\documentclass[11pt,twoside,reqno]{amsart}
\usepackage{srcltx} % for inverse search in LINUX -- does not cause any problems in Windows
\usepackage{amsmath,amsfonts,amsthm,amssymb}
\usepackage{graphicx,psfrag,overpic}
\usepackage{bm}
\usepackage[all,cmtip,line]{xy}
\usepackage{array}
\usepackage{enumitem}
\newtheorem{problemL}{Problem}
\newtheorem{theorem}{Theorem}%[section]
\newtheorem{lemma}{Lemma}[section]
\newtheorem{remark}[lemma]{Remark}
\newtheorem{proposition}[lemma]{Proposition}
\theoremstyle{definition}

\newcommand \gam{\gamma}
\newcommand \R{\mathbb{R}}
\newcommand \Om{\Omega}
\newcommand \der{\partial}
\newcommand \vphi{\varphi}
\newcommand \mcl{\mathcal}
\newcommand \Gam{\Gamma}
\newcommand \alp{\alpha}

\newcommand \til{\tilde}
\newcommand \wtil{\widetilde}
\newcommand \ol{\overline}

\newcommand \eps{\varepsilon}

\newcommand{\x}{\mathbf{x}}
\newcommand{\nnu}{{\boldsymbol\nu}}

\newcommand \Shock{\Gamma_{shock}}
\newcommand \Wedge{\Gamma_{wedge}}

\newcommand \shock{\Gamma_{shock}}

\newcommand \lefttop{P_1}
\newcommand \righttop{P_2}

\newcommand \leftbottom{P_4}
\newcommand \rightbottom{P_3}

\newcommand \leftshock{S_{\mathcal{O}}}
\newcommand \rightshock{S_{\mcl{N}}}
\newcommand \leftsonic{\Gamma_{sonic}^{\mathcal{O}}}
\newcommand \rightsonic{\Gamma_{sonic}^{\mcl{N}}}
\newcommand \leftvec{\hat{e}_{S_{\mathcal{O}}}}
\newcommand \rightvec{-\hat{\xi}}
\newcommand \leftvphi{\varphi_{\mathcal{O}}}
\newcommand \rightvphi{\varphi_{\mcl{N}}}
\newcommand \ivphi{\varphi_{\infty}}
\newcommand \iv{v_{\infty}}
\newcommand \iu{u_{\infty}}

\newcommand \irho{\rho_{\infty}}

\newcommand \leftu{u_{\mathcal{O}}} \newcommand \rightu{u_{\mathcal{N}}}
\newcommand \leftv{v_{\mathcal{O}}} \newcommand \rightv{v_{\mathcal{N}}}
\newcommand \leftc{c_{\mathcal{O}}}
\newcommand \leftrho{\rho_{\mathcal{O}}}
\newcommand \rightc{c_{\mcl{N}}}
\newcommand \rightrho{\rho_{\mcl{N}}}
\newcommand \ik{k_{\infty}}
\newcommand \leftk{k_{\mcl{O}}}
\newcommand \rightk{k_{\mcl{N}}}
\newcommand \nrho{\rho_{\mcl{N}}}

\newcommand \oxi{\xi_{\mcl{O}}}
\newcommand \nxi{\xi_{\mcl{N}}}

\newcommand \neta{\eta_{\mcl{N}}}
\newcommand \oeta{\eta_{\mcl{O}}}

\newcommand \oL{L_{\mcl{O}}}
\newcommand \iL{L_{\infty}}

\newcommand \otheta{\theta_{\mcl{O}}}
\newcommand \itheta{\theta_{\infty}}

\newcommand \nD{\mcl{D}^{\mcl{N}}}
\newcommand \oD{\mcl{D}^{\mcl{O}}}

\newcommand \oOm{\Omega^{\mcl{O}}}

\newcommand \betac{\beta_{sonic}}

\newcommand \cone{\mathrm{Cone}^0(\leftvec, \rightvec)}

\newcommand \iter{\mcl{Q}^{iter}}

\newcommand \fshock{f_{shock}}

\numberwithin{equation}{section}

\begin{document}

\title[]
{Prandtl-Meyer Reflection \\for Supersonic Flow past a Solid Ramp}
\author{Myoungjean Bae}
\address{M. Bae, Department of Mathematics, POSTECH,
         San 31, Hyojadong, NAmgu, Pohang, Gyungbuk, Korea}
\email{mjbae@postech.ac.kr}
\author{Gui-Qiang Chen}
\address{G.-Q. Chen, Mathematical Institute, University of Oxford, 24-29 St Giles', Oxford, OX1 3LB, England;
School of Mathematical Sciences, Fudan University, Shanghai 200433, China;
Department of Mathematics, Northwestern University,
2033 Sheridan Road, Evanston, IL 60208-2734, USA}
\email{chengq@maths.ox.ac.uk}

\author{Mikhail Feldman}
\address{M. Feldman, Department of Mathematics, University of Wisconsin, Madison, WI 53706-1388, USA}
\email{feldman@math.wisc.edu}

\bigskip
\dedicatory{Dedicated to Costas Dafermos on the occasion of his 70th
birthday}

\date{\today}
\maketitle

\begin{abstract}
We present our recent results on the Prandtl-Meyer reflection for supersonic potential
flow past a solid ramp.
When a steady supersonic flow passes a solid ramp,
there are two possible configurations:
the weak shock solution and the strong shock solution.
Elling-Liu's theorem (2008) indicates that the steady supersonic weak shock solution
can be regarded as a long-time asymptotics of an unsteady flow
for a class of physical parameters determined by certain
assumptions for potential flow.
In this paper we discuss our recent progress in removing these assumptions and establishing
the stability theorem for steady supersonic weak shock solutions as
the long-time asymptotics of unsteady flows
for all the physical parameters for potential flow.
We apply new mathematical techniques developed in our recent work to obtain
monotonicity properties and uniform apriori estimates for weak solutions,
which allow us to employ
the Leray-Schauder degree argument to complete the theory for the general case.
\end{abstract}

\medskip
\section{Introduction}

We are concerned with unsteady global solutions for supersonic flow past a solid ramp,
which may be also regarded as portraying the symmetric gas flow impinging onto a solid wedge (by symmetry).
When a steady supersonic flow past a solid ramp whose slope is less than a critical slope,
Prandtl employed the shock polar analysis to show that there are two possible configurations:
the weak shock reflection with supersonic downstream flow  and the strong shock reflection
with subsonic downstream flow, which both satisfy the entropy
conditions, provided that we do not give additional conditions at downstream;
see Busemann \cite{Busemann}, Courant-Friedrichs \cite{CF}, Meyer \cite{Meyer},
and Prandtl \cite{Prandtl}.

The fundamental question of whether one or both of the strong and the weak shocks are physically
admissible has been vigorously debated over the past seventy years, but has not yet been settled
in a definite manner (cf. Courant-Friedrichs \cite{CF}, Dafermos \cite{Dafermos}, and Serre \cite{Serre}).
On the basis of experimental and numerical evidence, there are strong indications that it is
the Prandtl-Meyer weak reflection solution that is physically admissible.
One plausible approach is to single out the strong shock reflection by the consideration
of stability, the stable ones are physical.
It has been shown in the steady regime that the Prandtl-Meyer weak reflection
is not only structurally stable (cf. Chen-Zhang-Zhu \cite{CZZ}),
but also $L^1$-stable with respect to steady small perturbation of both the ramp slope
and the incoming steady upstream flow (cf. Chen-Li \cite{CLi}),
while the strong reflection is also structurally stable for a large spectrum of
physical parameters (cf. Chen-Fang \cite{ChenFang}).
The first rigorous unsteady analysis of the steady supersonic weak shock solution
as the long-time behavior of an unsteady flow
is due to Elling-Liu in \cite{EL2}
in which they
succeeded to establish a stability theorem for a class of physical
parameters determined by certain assumptions for potential flow (see \S 3.1).

The purpose of this work is to remove the assumptions in Elling-Liu's theorem \cite{EL2}
and establish the stability theorem for the steady supersonic weak shock solution as
the long-time asymptotics of an unsteady flow
for all the admissible physical parameters (without additional conditions)
for potential flow.
To complete this theorem, we apply new mathematical techniques developed in Chen-Feldman
\cite{CF2} to obtain uniform
apriori estimates for weak solutions.
We first establish the monotonicity property and its consequence of weak solutions (see \S 4.1).
Then we make various uniform apriori estimates of weak solutions for two cases:
the $C^{2,\alpha}$--estimates away from the sonic circles where the governing equation
is uniformly elliptic (see \S 4.3.1) and the weighted $C^{2,\alpha}$--estimates near the sonic
circles where the ellipticity degenerates (see \S 4.3.2).
These careful estimates allow us to employ
the Leray-Schauder degree argument to establish the complete theory (see \S 4.4).

\section{Mathematical Formulation for the Problem}
\label{section:shock-polar}

The unsteady potential flow is governed by the conservation law of mass and Bernoulli's law:
\begin{align}
\label{1-a}
&\der_t\rho+\nabla_{\bf x}\cdot(\rho \nabla_{\bf x}\Phi)=0,\\
\label{1-b}
&\der_t\Phi+\frac 12|\nabla_{\bf x}\Phi|^2+i(\rho)=B
\end{align}
for the density $\rho$ and the velocity potential $\Phi$, where the Bernoulli constant
$B$ is determined by the incoming flow and/or boundary conditions,
and $i(\rho)$ satisfies the relation
\begin{equation*}
i'(\rho)=\frac{p'(\rho)}{\rho}=\frac{c^2(\rho)}{\rho}
\end{equation*}
with $c(\rho)$ being the sound speed, and $p$ is the pressure that is a function of the density $\rho$.
For an ideal polytropic gas, the pressure $p$ and the sound speed $c$ are given by
\begin{equation*}
p(\rho)=\kappa \rho^{\gam}, \quad c^2(\rho)=\kappa \gam \rho^{\gam-1}
\end{equation*}
for constants $\gam>1$ and $\kappa>0$.
Without loss of generality, we choose $\kappa=1/\gam$ to have
\begin{equation}
\label{1-c}
i(\rho)=\frac{\rho^{\gam-1}-1}{\gam-1},\quad c^2(\rho)=\rho^{\gam-1}.
\end{equation}
This can be achieved by the following scaling:
\begin{equation*}
({\bf x}, t, B)\longrightarrow (\alp{\bf x}, \alp^2t, \alp^{-2}B),\quad \alp^2=\kappa\gam.
\end{equation*}
Taking the limit $\gam\to 1+$, we can also consider the case of the isothermal flow ($\gam=1$).
For isothermal flow, \eqref{1-c} implies
\begin{equation*}
i(\rho)=\ln \rho,\qquad c^2(\rho)\equiv 1.
\end{equation*}

Our goal is to find a solution $(\rho, \Phi)$ to system \eqref{1-a}--\eqref{1-b} when,
at $t=0$, a uniform flow in $\R^2_+:=\{x_1\in \R, x_2>0\}$ with
$(\rho, \nabla_{\bf x}\Phi)=(\irho, \iu,0)$ is heading to a solid ramp:
$$
W:=\{\x=(x_1,x_2): 0<x_2<x_1\tan\theta_w, x_1>0\}.
$$

\medskip
\begin{problemL}
[Initial-Boundary Value Problem]
\label{problem-1}
Seek a solution of system \eqref{1-a}--\eqref{1-b} with
$B=\frac{\iu^2}{2}+\frac{\irho^{\gam-1}-1}{\gam-1}$ and the initial condition at $t=0$:
\begin{equation}
\label{1-d}
(\rho,\Phi)|_{t=0}=(\irho, \iu x_1)
\qquad\text{for}\;\;(x_1,x_2)\in \R^2_+\setminus W,
\end{equation}
and with the slip boundary condition along the wedge boundary $\der W$:
\begin{equation}
\label{1-e}
\nabla_{\bf x}\Phi\cdot \nnu|_{\der W\cap\{x_2>0\}}=0,
\end{equation}
where $\nnu$ is the exterior unit normal to $\der W$.
\end{problemL}

\medskip
Notice that the initial-boundary value problem \eqref{1-a}--\eqref{1-b}
with \eqref{1-d}--\eqref{1-e} is invariant under the scaling:
\begin{equation*}
({\bf x}, t)\rightarrow (\alp{\bf x}, \alp t),\quad (\rho, \Phi)\rightarrow (\rho, \frac{\Phi}{\alp})\qquad \text{for}\;\;\alp\neq 0.
\end{equation*}
Thus, we seek self-similar solutions in the form of
\begin{equation*}
\rho({\bf x},t)=\rho(\xi,\eta),\quad \Phi({\bf x},t)=t\phi(\xi,\eta)\qquad\text{for}\;\;(\xi,\eta)=\frac{{\bf x}}{t}.
\end{equation*}
Then the pseudo-potential function $\vphi=\phi-\frac 12(\xi^2+\eta^2)$ satisfies
the following Euler equations for self-similar solutions:
\begin{align}
\label{1-r}
&{\rm div}(\rho D\vphi)+2\rho=0,\\
\label{1-p1}
&\frac{\rho^{\gam-1}-1}{\gam-1}+(\frac 12|D\vphi|^2+\vphi)=B,
\end{align}
where the divergence ${\rm div}$ and gradient $D$ are with respect to $(\xi,\eta)$.
From this, we obtain an equation for the pseudo-potential function $\vphi(\xi,\eta)$ as follows:
\begin{equation}
\label{2-1}
{\rm div}(\rho(|D\vphi|^2,\vphi)D\vphi)+2\rho(|D\vphi|^2,\vphi)=0
\end{equation}
for
\begin{equation}
\label{1-o}
\rho(|D\vphi|^2,\vphi)=
\bigl(B_0-(\gam-1)(\frac 12|D\vphi|^2+\vphi)\bigr)^{\frac{1}{\gam-1}},
\end{equation}
where we set $B_0:=(\gam-1)B+1$.
Then we have
\begin{equation}
\label{1-a1}
c^2(|D\vphi|^2,\vphi)=
B_0-(\gam-1)(\frac 12|D\vphi|^2+\vphi).
\end{equation}
Equation \eqref{2-1} is an equation of mixed elliptic-hyperbolic type. It is elliptic if and only if
\begin{equation}
\label{1-f}
|D\vphi|<c(|D\vphi|^2,\vphi).
\end{equation}

If $\rho$ is a constant, then, by \eqref{2-1} and \eqref{1-o}, the corresponding
pseudo-potential $\vphi$ is in the form of
$$
\vphi(\xi,\eta)=-\frac 12(\xi^2+\eta^2)+u\xi+v\eta+k
$$
for constants $u,v$, and $k$.
As the incoming flow has the constant velocity $(\iu,0)$, the corresponding pseudo-potential
$\ivphi$ has the expression of
\begin{equation}
\label{1-m}
\ivphi=-\frac 12(\xi^2+\eta^2)+\iu\xi+\ik
\end{equation}
for a constant $\ik$.

Without loss of generality by scaling, we fix $\irho=1$ and $M_\infty=\frac{\iu}{\irho^{(\gamma-1)/2}}=\iu$ (Mach
number of the state at infinity), and then \eqref{1-p1} becomes
\begin{equation}
\label{1-p}
\frac{{\rho}^{\gam-1}-1}{\gam-1}+\frac 12(|D{\vphi}|^2+{\vphi})=\frac{{M_\infty}^2}{2}.
\end{equation}
Then Problem \ref{problem-1} can be reformulated as the following boundary value problem in the
self-similar coordinates $(\xi,\eta)$.  The domain in the self-similar coordinates $(\xi,\eta)$
corresponding to $\{({\bf x}, t)\, :\, {\bf x}\in \R^2_+\setminus W,\, t>0\}$ is
$$
\Lambda:=\R^2_+\setminus\{(\xi,\eta):\eta\le \xi\tan\theta_w,\, \xi\ge 0\}.
$$

\medskip
\begin{problemL}[Boundary Value Problem]
\label{problem-2}
Seek a solution $\vphi$ of equation \eqref{2-1} in the self-similar domain $\Lambda$
with the slip boundary condition:
\begin{equation}
\label{1-k}
D\vphi\cdot\nnu_w=0\qquad\text{on}\;\;\Wedge=\{(\xi,\eta): \eta=\xi\tan\theta_w, \xi>0\},
\end{equation}
where $\nnu_w$ is the exterior unit normal to the boundary of the wedge $\Wedge$.
\end{problemL}

\begin{figure}[htp]
\centering
\begin{psfrags}
\psfrag{C}[cc][][0.8][0]{$c_*$}
\psfrag{TC}[cc][][0.8][0]{$\theta_c$}
\psfrag{TW}[cc][][0.8][0]{$\theta_w$}
\psfrag{U}[cc][][0.8][0]{$u$}
\psfrag{V}[cc][][0.8][0]{$v$}
\psfrag{U0}[cc][][0.8][0]{$u_0$}
\psfrag{I}[cc][][0.8][0]{$\Om_{\infty}:(\iu, 0), p_{\infty}, \rho_{\infty}$}
\psfrag{TH}[cc][][0.8][0]{$\theta_w$}
\psfrag{Om}[cc][][0.8][0]{$\Omega$}
\psfrag{L}[cc][][0.8][0]{$\Om_{\mcl{O}}$}
\psfrag{R}[cc][][0.8][0]{$\Om_{\mcl{N}}$}
\psfrag{S2}[cc][][0.8][0]{$\leftshock$}
\psfrag{S1}[cc][][0.8][0]{$\rightshock$}
\psfrag{S}[cc][][0.8][0]{$S$}
\includegraphics[scale=.4]{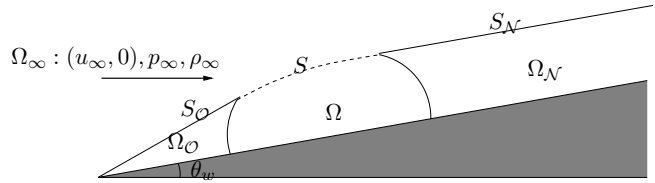}
\caption{Weak shock solutions in the self-similar coordinates}\label{fig:global-structure}
\end{psfrags}
\end{figure}

In particular, we seek a weak solution of Problem \ref{problem-2} such that it contains
a straight weak oblique shock attached to the tip of the wedge, and the oblique shock is
connected to a normal shock through a curved shock as shown in Figure \ref{fig:global-structure}.

A shock is a curve across which $D\vphi$ is discontinuous. If $\Om^+$ and $\Om^-(:=\Om\setminus \ol{\Om^+})$
are two nonempty open subsets of $\Om\subset \R^2$, and $S:=\der\Om^+\cap \Om$ is a $C^1$-curve
where $D\vphi$ has a jump, then $\vphi\in W^{1,1}_{loc}\cap C^1(\Om^{\pm}\cup S)\cap C^2(\Om^{\pm})$
is a global weak solution of \eqref{2-1} in $\Om$ if and only if $\vphi$ is in $W^{1,\infty}_{loc}(\Om)$
and satisfies equation \eqref{2-1} and the Rankine-Hugoniot condition on $S$:
\begin{equation}
\label{1-h}
[\rho(|D\vphi|^2, \vphi)D\vphi\cdot\nnu]_S=0,
\end{equation}
where $[F]_S$ is defined by
$$
[F(\xi,\eta)]_S:=F(\xi,\eta)|_{\overline{\Om^-}}-F(\xi, \eta)|_{\overline{\Om^+}}\quad\text{for}\;\;(\xi, \eta)\in S.
$$
Note that the condition $\vphi\in W^{1,\infty}_{loc}(\Om)$ requires
\begin{equation}
\label{1-i}
[\vphi]_S=0.
\end{equation}
In Figure \ref{fig:global-structure}, by \eqref{1-m} and \eqref{1-i},
the pseudo-potentials $\leftvphi$ and $\rightvphi$ below $\leftshock$ and $\rightshock$ are respectively
in the form of
\begin{equation}
\label{1-n}
\begin{split}
&\leftvphi=-\frac 12(\xi^2+\eta^2)+\leftu \xi+\leftv\eta+\leftk,\\
&\rightvphi=-\frac 12(\xi^2+\eta^2)+\rightu \xi+\rightv \eta+\rightk
\end{split}
\end{equation}
for constants $\leftu, \leftv, \rightu, \rightv, \leftk$, and $\rightk$. Then it follows
from \eqref{1-o} and \eqref{1-n} that the corresponding densities $\leftrho$ and $\rightrho$
below $\leftshock$ and $\rightshock$ are constants, respectively.

\begin{figure}[htp]
\centering
\begin{psfrags}
\psfrag{tc}[cc][][0.8][0]{$\theta_c$}
\psfrag{tw}[cc][][0.8][0]{$\theta_w$}
\psfrag{u}[cc][][0.8][0]{$u$}
\psfrag{v}[cc][][0.8][0]{$v$}
\psfrag{u0}[cc][][0.8][0]{$\iu$}
\psfrag{zeta}[cc][][0.8][0]{}
\includegraphics[scale=.6]{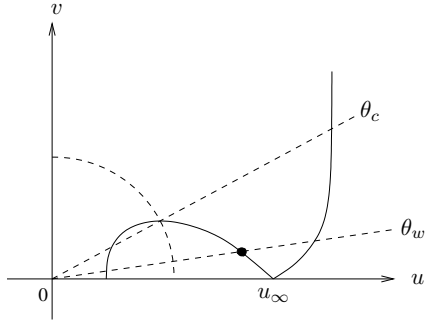}
\caption{Shock polar}\label{fig:polar}
\end{psfrags}
\end{figure}

Given $M_\infty>1$, we obtain $(\leftu, \leftv)$ and $\leftrho$ by using the shock polar curve
for steady potential flow in Figure \ref{fig:polar}. In Figure \ref{fig:polar},
let $\theta_{sonic}$ be the wedge angle such that the line $\frac{v}{u}=\tan\theta_{sonic}$ intersects
with the shock polar curve at a point on the circle of radius $1$.
For a wedge angle $\theta_w\in(0,\theta_{sonic})$, the line $v=u\tan\theta_w$ and the shock polar
curve intersect at a point $(\leftu, \leftv)$ with $\sqrt{\leftu^2+\leftv^2}>1$ and $\leftu<\iu$.
The intersection $(\leftu, \leftv)$ indicates the velocity for steady potential flow
behind an oblique shock $\leftshock$ attached to the tip of the wedge angle $\theta_w$.
Since the strength of the shock $\leftshock$ is relatively weak compared to the other shock
given by the other intersection point on the shock polar curve,
$\leftshock$ is called a \emph{weak shock}.
We also note that the states on both sides of $\leftshock$ are supersonic,
and such states $(\leftu, \leftv)$ smoothly depend on $\iu$ and $\theta_w$.
Once $(\leftu, \leftv)$ is determined, then, by \eqref{1-p},
the corresponding density $\leftrho$ is given by
\begin{equation}
\label{2-n1}
\leftrho^{\gam-1}=1+\frac{\gam-1}{2}(\iu^2-\leftu^2-\leftv^2).
\end{equation}

To find a curved shock $S$ connecting the oblique shock $\leftshock$ with the normal shock $\rightshock$
in Figure \ref{fig:global-structure}, it is more convenient to use a coordinate system in which
the boundary of wedge $\Wedge$ is the same for all $\theta_w\in(0,\theta_{sonic})$.
For that reason, we introduce a new coordinate system $(\xi',\eta')$ defined by
\begin{equation}
\label{2-a}
\begin{pmatrix}
\xi'\\
\eta'
\end{pmatrix}
=\begin{pmatrix}
\cos\theta_w&\sin\theta_w\\
-\sin\theta_w&\cos\theta_w
\end{pmatrix}
\begin{pmatrix}
\xi\\
\eta
\end{pmatrix}-
\begin{pmatrix}
\leftu\cos\theta_w\\
0\end{pmatrix}.
\end{equation}
In other words, $(\xi',\eta')$ is obtained by rotating $(\xi, \eta)$ and
translating the self-similar plane.
Denoting as $(\xi,\eta)$ for $(\xi',\eta')$,
$\ivphi$ and $\leftvphi$ defined by \eqref{1-m} and \eqref{1-n} are expressed as
\begin{equation}
\label{1-21}
\begin{split}
&\ivphi=-\frac 12(\xi^2+\eta^2)-\eta \iu\sin\theta_w+\ik,\\
&\leftvphi=
-\frac 12(\xi^2+\eta^2)+\wtil{\leftu}\xi+\leftk
\end{split}
\end{equation}
with
\begin{equation}
\label{J}
\wtil{\leftu}=\sqrt{\leftu^2+\leftv^2}-\iu\cos\theta_w
\end{equation}
for constants $\ik$ and $\leftk$ different from \eqref{1-m} and \eqref{1-n}.
For simplicity, we will write $\wtil{\leftu}$ as $\leftu$ hereafter.
Without loss of generality, we choose $\ik=0$.
Set
\begin{equation}
\label{1-22}
\iv:=\iu \sin\theta_w,
\end{equation}
and let $\beta$ be the angle between the oblique shock $\leftshock$ and wedge $\Wedge$.
We use the parameters $(\iv,\beta)$ instead of $(\iu, \theta_w)$ to compute
the normal shock $\rightshock$ and the state behind it.

\begin{figure}[htp]
\centering
\begin{psfrags}
\psfrag{iv}[cc][][0.7][0]{$(0,-\iv)$}
\psfrag{eta}[cc][][0.7][0]{$\eta$}
\psfrag{nshock}[cc][][0.7][0]{$\rightshock$}
\psfrag{normal}[cc][][0.7][0]{$\eta=\neta$}
\psfrag{lshock}[cc][][0.7][0]{$\leftshock$}
\psfrag{shock}[cc][][0.7][0]{$\Shock$}
\psfrag{lt}[cc][][0.7][0]{$\lefttop$}
\psfrag{rt}[cc][][0.7][0]{$\righttop$}
\psfrag{osonic}[cc][][0.7][0]{$\leftsonic$}
\psfrag{nsonic}[cc][][0.7][0]{$\rightsonic$}
\psfrag{B}[cc][][0.7][0]{$\beta$}
\psfrag{tip}[cc][][0.7][0]{$\til{\xi}(\beta)$}
\psfrag{lb}[cc][][0.7][0]{$\leftbottom$}
\psfrag{rb}[cc][][0.7][0]{$\rightbottom$}
\psfrag{O}[cc][][0.7][0]{0}
\psfrag{xi}[cc][][0.7][0]{$\xi$}
\includegraphics[scale=.4]{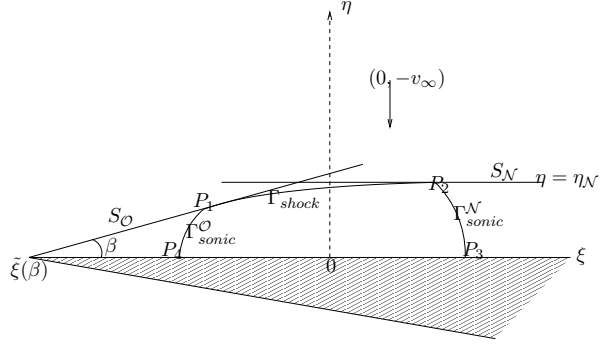}
\caption{Weak shock solutions on the self-similar plane}\label{fig:regularsol1}
\end{psfrags}
\end{figure}

By \eqref{1-i}, the definition of $\ivphi$ given by \eqref{1-21} with $\ik=0$,
and by the slip boundary condition \eqref{1-k},
if we let $\rightshock=\{\eta=\neta\}$ for a constant $\neta$, then we have
\begin{equation}
\label{2-k}
\rightvphi=-\frac 12(\xi^2+\eta^2)-\iv\neta.
\end{equation}
By \eqref{1-h} and \eqref{1-p}, the density $\rightrho$ and constant $\neta$ satisfy
\begin{align}
\label{1-23}
&\neta=\frac{\iv}{\rightrho-1},\\
\label{1-1a}
&\frac{\rightrho^{\gam-1}-1}{\gam-1}=\frac 12\iv^2+\neta \iv.
\end{align}
We note that $0<\neta<\rightc$ for $\rightc=\nrho^{(\gam-1)/2}$.

Fix $\beta\in(0,\frac{\pi}{2})$, consider an oblique shock $\leftshock$ of the angle $\beta$
from the $\xi-$axis, and let $\leftvphi$ be the corresponding pseudo-potential below $\leftshock$.
Let $(\oxi,\oeta)$ be the intersection $\lefttop$ of $\leftshock$ and
the sonic circle $B_{\leftc}(\leftu,0)$ for the state in $\oOm$ (see Figure 1).
There are two intersection points, and we denote $\lefttop$ the point with
the smaller value of the $\xi-$coordinate (see Figure \ref{fig:regularsol1}).
Given $\iv>0$ and $\gam\ge 1$,  $\oeta$  satisfies
\begin{equation}
\label{1-b1}
\frac{\der\oeta}{\der\beta}<0\qquad\text{for all}\;\;\beta\in (0,\frac{\pi}{2}).
\end{equation}
 Therefore, the set
$
I_{\iv}:=\{\beta\in(0,\frac{\pi}{2}):0<\oeta(\beta)<\infty\}
$
is connected and, furthermore, there exists $\betac\in(0,\frac{\pi}{2}]$ satisfying
$$
I_{\iv}=(0,\betac),
$$
and such $\betac$ depends smoothly on $\iv$ and $\gam$.

Given  $\gam\ge 1$, we note that $\theta_{sonic}$ on the shock polar curve in Figure \ref{fig:polar}
depends on $\iu$, while $\betac$ depends on $\iv$, so that we write as $\theta_{sonic}(\iu)$
and $\betac(\iv)$, respectively.
Define two parameter sets $\mathfrak{P}_1$ and $\mathfrak{P}_2$ by
\begin{equation}
\label{1-38}
\begin{split}
&\mathfrak{P}_1:=\bigcup_{\iu>1} \{\iu\}\times (0,\theta_{sonic}(\iu)),\\
&\mathfrak{P}_2:=\bigcup_{\iv>0} \{\iv\}\times (0,\betac(\iv)).
\end{split}
\end{equation}
Then we have the following lemma.
\begin{lemma}
\label{lemma-parameters}
For any $\gam\ge 1$,
there exists a one-to-one and onto correspondence between $\mathfrak{P}_1$ and $\mathfrak{P}_2$.
\end{lemma}

\begin{remark}[Normal shock: $\beta=0$]
For fixed $\gam\ge 1$ and $\iv>0$, the state of $\beta=0$ is a normal shock,
and this corresponds to the state of $(\iu, \theta_w)=(\infty,0)$.
Even though $\beta=0$ is nonphysical, we will put the case of $\beta=0$ in our consideration
as it is useful in applying the Leray-Schauder degree theorem to establish
Theorem {\rm \ref{theorem-1}} stated in Section {\rm \ref{subsec-mainthm}}.
\end{remark}

For a fixed  $(\iv,\beta)\in \mathfrak{P}_2$,  if we can prove the existence of
a curved shock $\shock$ to connect $\leftshock$ with $\rightshock$
and the existence of corresponding pseudo-potential below $\shock$,
then, by Lemma \ref{lemma-parameters}, the existence of a curved shock and
the pseudo-potential below the shock for $(\iu,\theta_w)\in\mathfrak{P}_1$
corresponding to the fixed $(\iv, \beta)$ automatically follows.
By taking $(\iv,\beta)$ as a parameter instead of $(\iu,\theta_w)$,
it is advantageous that the boundary of a wedge is always on the $\xi$-axis.
As $(\iv, \beta)$ is the parameter of our problem,
we write $\ivphi$ in \eqref{1-21} as
\begin{equation}
\label{2-a1}
\ivphi=-\frac 12(\xi^2+\eta^2)-\iv \eta.
\end{equation}

Since the states below $\leftshock$ and $\rightshock$ are given by $\leftvphi$ and $\rightvphi$,
in order to solve Problem \ref{problem-2}, it suffices to solve the following problem.

\medskip
\begin{problemL}[Free Boundary Problem]\label{fbp}
Find a curved shock $\shock$ and a function $\vphi$ defined in the region $\Om$, enclosed
by $\shock, \leftsonic, \rightsonic$, and $\{\eta=0\}$, such that $\vphi$ satisfies
\begin{itemize}
\item[\rm (i)]
Equation \eqref{2-1} in $\Om$;
\item[\rm (ii)]
$\vphi=\ivphi$, $\rho D\vphi\cdot\nnu_s=D\ivphi\cdot\nnu_s$ {on} $\shock$;
\item[\rm (iii)]
$\vphi=\vphi_{\beta}$, $D\vphi=D\vphi_{\beta}$ {on} $\leftsonic\cup\rightsonic\,\,$
 for $\vphi_{\beta}:=\max(\leftvphi, \rightvphi)$;
\item[\rm (iv)]
$\der_{\eta}\vphi=0$ {on} $\Wedge$,
\end{itemize}
where $\nnu_s$ is the interior unit normal on $\shock$.
\end{problemL}

Let $\vphi$ be a solution of Problem \ref{fbp} with a shock $\shock$.
Moreover, assume that $\vphi\in C^1(\overline\Om)$, and $\shock$ is a $C^1$--curve up to its endpoints.
To obtain a solution of Problem \ref{problem-2} from $\vphi$,
we divide the half-plane $\{\eta\ge 0\}$ into four separate regions.
Let $\Om_{E}$ be the unbounded domain below the curve $\ol{\leftshock\cup\shock\cup\rightshock}$
in $\{\eta\ge 0\}$ (see Figure \ref{fig:regularsol1}).
In $\Om_E$, let $\Om_{\mcl{O}}$ be the bounded open domain enclosed by $\leftshock, \leftsonic$,
and $\{\eta=0\}$.
We set $\Om_{\mcl{N}}:=\Om_E\setminus \ol{(\Om_{\mcl{O}}\cup\Om)}$.
Define a function $\vphi_*$ in $\{\eta\ge 0\}$ by
\begin{equation}\label{extsol}
\vphi_*=
\begin{cases}
\ivphi& \qquad \text{in}\,\{\eta\ge 0\}\setminus \Om_E,\\
\leftvphi& \qquad \text{in}\,\Om_{\mcl{O}},\\
\vphi& \qquad \text{in}\, \leftsonic\cup\Om\cup\rightsonic,\\
\rightvphi&\qquad \text{in}\,\Om_{\mcl{N}}.
\end{cases}
\end{equation}
By \eqref{1-i} and (iii) of Problem \ref{fbp}, $\vphi_*$ is continuous in $\{\eta\ge 0\}$
and is $C^1$ in $\overline{\Om_E}$.
In particular, $\vphi_*$ is $C^1$ across $\leftsonic\cup\rightsonic$.
Moreover, using (i)--(iii) of Problem \ref{fbp}, we obtain
that  $\vphi_*$ is a weak solution of equation \eqref{2-1}
in $\{\eta>0\}$.
Applying the inverse coordinate transformation of \eqref{2-a}
to $\vphi_*$, we obtain a solution of Problem \ref{problem-2}.

\section{Results}\label{subsec-mainthm}

In this section we first present the known result, Elling-Liu's theorem, in \cite{EL2}, and then we describe
our new results for Problem 3, the Prandtl-Meyer problem.

\subsection{Known Result}

Elling and Liu in \cite{EL2} have established the following theorem.

\medskip
\noindent
{\bf Elling-Liu's Theorem}
(Theorem 1 in \cite{EL2}).
Given $(\iv,\beta)\in\mathfrak{P}_2$, let $L$ be the line segment connecting
$\lefttop$ with $\righttop$ in Figure \ref{fig:regularsol1}.
If $L$ does not intersect with $B_{1}(0,-\iv)$, then there exists a weak
solution to Problem \ref{problem-2}
with structure \eqref{extsol}.

\begin{figure}[htp]
\centering
\begin{psfrags}
\psfrag{ti}[cc][][0.7][0]{$\itheta$}
\psfrag{to}[cc][][0.7][0]{$\otheta$}
\psfrag{n}[cc][][0.7][0]{$(\nxi,\neta)$}
\psfrag{o}[cc][][0.7][0]{$(\oxi,\oeta)$}
\psfrag{sn}[cc][][0.7][0]{$\rightshock$}
\psfrag{so}[cc][][0.7][0]{$\leftshock$}
\psfrag{iv}[cc][][0.7][0]{$\;\;-\iv$}
\psfrag{d0}[cc][][0.7][0]{$d=1\;\;\;$}
\includegraphics[scale=.6]{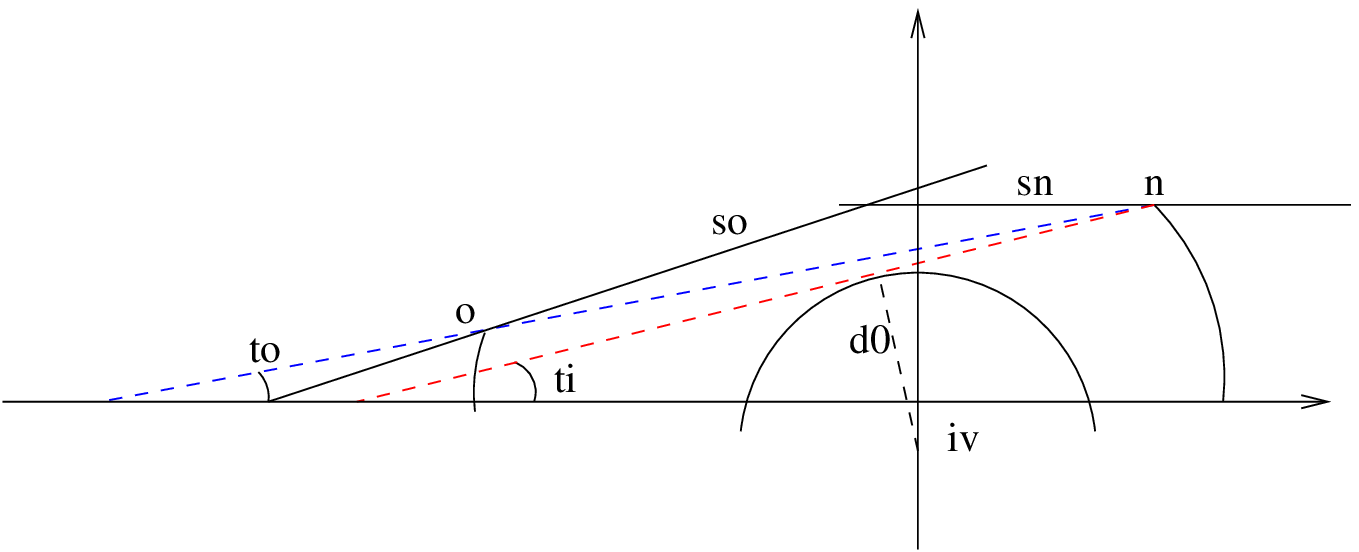}
\includegraphics[scale=.6]{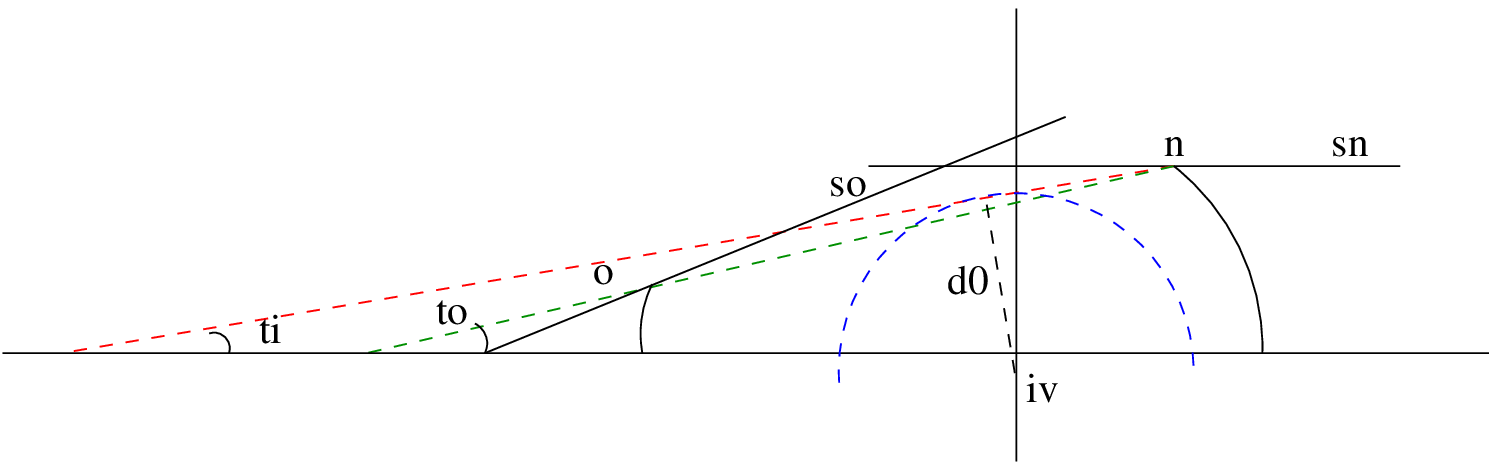}
\caption{Top: $\otheta<\itheta$; $\quad$ Bottom: $\otheta>\itheta$ }
\label{fig:el}
\end{psfrags}
\end{figure}

\medskip
If $\iv>1$, since $B_1(0,-\iv)\subset\{\eta<0\}$, the assumption
in the theorem above holds true.
On the other hand, for $\iv<1$, there exists a set of parameters  $(\iv, \beta)\in \mathfrak{P}_2$
for which the assumption no long holds. This can be shown as follows:

Fix $\gam\ge 1$ and $\iv>0$.
Let $\iL$ be the line with positive slope such that $\iL$ is tangent to
the sonic circle $B_{1}(0,-\iv)$ of the incoming state and
that it passes through the point $\righttop=(\nxi,\neta)=\ol{\rightsonic}\cap \rightshock$.
For some $\iv>0$, there may be two such tangent lines.
In that case, we choose the one with the smaller slope.
We also fix $\beta\in(0,\betac(\iv))$.
Let $\oL$ be the line segment
connecting $\lefttop:=(\oxi, \oeta)=\ol{\leftsonic}\cap \leftshock$
with $\righttop$ (see Figure \ref{fig:el}).

Let $\tan\otheta$ and $\tan\itheta$ be the slopes of the lines $\oL$ and $\iL$,
respectively. If $\otheta>\itheta$, then the line $\oL$ intersects with the sonic
circle $B_{1}(0,-\iv)$. On the other hand, if $\otheta<\itheta$,
then $\oL$ has no intersection with $B_{1}(0,-\iv)$. Define
\begin{equation}
\label{1-17}
F(\beta):=\tan\otheta-\tan\itheta.
\end{equation}
Then, we have the following proposition.

\begin{proposition}
\label{proposition:el-condition}
For any given $\gamma\ge 1$, there exists $v_*>0$ depending on $\gam$ so that,
wherever $\iv\in(0,v_*)$, there exists
$\hat{\beta}=\hat{\beta}(\iv)\in(0,\betac(\iv))$ satisfying
\begin{equation}
\label{1-16}
\begin{split}
&F(\beta)\le0\qquad \text{for}\;\;\beta\in[0,\hat\beta],\\
&F(\beta)> 0\qquad \text{for}\;\;\beta\in(\hat\beta,\betac).
\end{split}
\end{equation}
\end{proposition}
According to Proposition \ref{proposition:el-condition},
there exists a subset of $\mathfrak{P}_2$ for which \cite[Theorem 1]{EL2}
does not apply.

\subsection{New Results}
Our new results provide not only the existence of a global solution
to Problem \ref{fbp} for all $(\iv, \beta)$ in $\mathfrak{P}_2$,
but also the higher regularity of the solution in pseudo-subsonic regions.
Therefore, we achieve the existence of a self-similar weak solution
with higher regularity to Problem \ref{problem-1}
for all $(\iu, \theta_w)$ in $\mathfrak{P}_1$.

\medskip
For any fixed $(\iv, \beta)\in\mathfrak{P}_2$, the oblique shock $\leftshock$ of
the slope $\tan\beta$, the normal shock $\rightshock$, and the corresponding
pseudo-potentials $\leftvphi$ and $\rightvphi$ below $\leftshock$ and $\rightshock$
are uniquely computed.
Let $(\til{\xi}(\beta),0)$ be the $\xi$-intercept of $\leftshock$.
Then we have
$$
\leftshock\subset\{(\xi,f_{\mcl{O}}(\xi)): f_{\mcl{O}}(\xi)=\tan\beta(\xi-\til{\xi}(\beta)), \xi \in \R\}.
$$
Denoting $\lefttop=(\oxi,\oeta)$ and $\righttop=(\nxi, \neta)$
in Figure \ref{fig:regularsol1}, we have
\begin{equation*}
\leftshock=\{(\xi,\eta): \til{\xi}(\beta)\le \xi\le \oxi, \eta=f_{\mcl{O}}(\xi)\},\quad\;
\rightshock=\{(\xi,\neta): \xi\ge\nxi\}.
\end{equation*}

\medskip
\begin{theorem}
\label{theorem-1}
For any given $\gam\ge 1$ and $(\iv,\beta)\in\mathfrak{P}_2$, there exists a global
weak solution $\vphi$ of {\rm Problem \ref{fbp}} satisfying the following properties:
\begin{itemize}
\item[\rm (i)] There exists a shock curve $\shock$ with endpoints $\lefttop=(\oxi,\oeta)$
and $\righttop=(\nxi,\neta)$ such that
\begin{itemize}
\item $\shock$ satisfies
$
\shock\subset(\R^2_+\setminus\ol{B_{1}(0,-\iv)})\cap\{\oxi\le \xi\le \nxi\};
$
\item $\shock$ is $C^3$ in its relative interior:
For any $P\in \shock\setminus\{\lefttop,\righttop\}$,
there exist $r>0$ and $f\in C^{3}(\R)$ such that
$$
\shock\cap B_r(P)=\{(\xi, \eta):\eta=f(\xi)\}\cap B_r(P);
$$
\item The curve $\overline{\leftshock\cup\shock\cup\rightshock}$ is $C^1$, including at the points $\lefttop$ and $\righttop$;
\item $\shock, \rightsonic, \leftsonic$, and $\Wedge:=\{(\xi,0):\leftu-\leftc\le\xi\le\rightc\}$
 do not have common points except at their end points.
 Thus, ${\shock}\cup {\rightsonic}\cup {\leftsonic}\cup{\Wedge}$ is a closed curve without
  self-intersections. Denote by $\Om$ the open and bounded domain enclosed by this closed curve.
\end{itemize}
\item[\rm (ii)] $\vphi$ satisfies $\vphi\in C^{3}(\ol{\Om}\setminus({\leftsonic}\cup{\rightsonic}))\cap C^1(\ol{\Om})$.
\item[\rm (iii)] Equation \eqref{2-1} is strictly elliptic in $\ol{\Om}\setminus({\rightsonic}\cup {\leftsonic})$.
\item[\rm (iv)] $\max(\leftvphi, \rightvphi)\le \vphi \le \ivphi$ in $\Om$.
\item[\rm (v)] $\der_{\rightvec}(\ivphi-\vphi)\le 0$ and $\der_{\leftvec}(\ivphi-\vphi)\le 0$ in $\Om$ for
\end{itemize}
\end{theorem}

If $\vphi$ is a weak solution satisfying properties (i)--(v) of Theorem \ref{theorem-1},
then $\vphi$ and $\shock$ satisfy additional regularity properties.

\begin{theorem} \label{theorem-3}
Given $\gam\ge 1$ and $(\iv, \beta)\in \mathfrak{P}_2$,
let $\vphi$ be a weak solution of {\rm Problem \ref{theorem-1}}
satisfying properties {\rm (i)}--{\rm (v)} of {\rm Theorem \ref{theorem-1}}.
Then the following properties hold:
\begin{itemize}
\item[\rm (i)] The curve $\overline{\leftshock\cup\shock\cup\rightshock}$ is $C^{2, \alpha}$
for any $\alpha\in [0, \frac{1}{2})$,
including at the points $\lefttop$ and $\righttop$.
Moreover, $\shock$ is $C^{\infty}$ in its relative interior.
\item[\rm (ii)] $\vphi\in C^{\infty}(\ol{\Om}\setminus({\leftsonic}\cup{\rightsonic}))\cap C^{1,1}(\ol{\Om})$.
\item[\rm (iii)]  For a constant $\sigma>0$ and a set $\mcl{D}$ given by
\begin{equation*}
\mcl{D}=\{ (\xi, \eta):\max(\leftvphi(\xi, \eta), \rightvphi(\xi, \eta))<\ivphi(\xi, \eta), \eta>0\},
\end{equation*}
define
\begin{equation*}
\begin{split}
&\nD_{\sigma}=\mcl{D}\cap \{(\xi,\eta): {\rm dist}((\xi,\eta),\rightsonic)<\sigma\}\cap B_{\rightc}(0,0),\\
&\oD_{\sigma}=\mcl{D}\cap \{(\xi,\eta): {\rm dist}((\xi,\eta),\leftsonic)<\sigma\}\cap B_{\leftc}(\leftu,0)\\
\end{split}
\end{equation*}
for $\rightc=\rightrho^{(\gam-1)/{2}}$ and $\leftc=\leftrho^{(\gam-1)/{2}}$.
Then, for any $\alp\in(0,1)$ and any given
$(\xi_0,\eta_0)\in (\ol{\leftsonic}\cup\ol{\rightsonic})\setminus\{\lefttop, \righttop\}$,
there exists $K<\infty$ depending only on
$\gam, \iv$, $\eps_0, \alp,$ $\|\vphi\|_{C^{1,1}(\Om\cap(\oD_{\eps_0}\cup\nD_{\eps_0}))}$,
and $d={\rm dist}((\xi_0,\eta_0),\shock)$ so that there holds
    \begin{equation}
    \label{Op1}
    \|\vphi\|_{2,\alp,\ol{\Om\cap B_{d/2}(\xi_0,\eta_0)\cap(\nD_{\eps_0/2}\cup\oD_{\eps_0/2})}}\le K.
    \end{equation}
\item[\rm (iv)] For any $(\xi_0,\eta_0)\in (\leftsonic\cup\rightsonic)\setminus\{\lefttop, \righttop\}$,
\begin{equation}
\label{Op2}
\lim_{(\xi,\eta)\to(\xi_0,\eta_0)\atop
(\xi,\eta)\in\Om}\big(D_{rr}\vphi-D_{rr}\max(\rightvphi,\leftvphi)\big)=\frac{1}{\gam+1},
\end{equation}
where $r=\sqrt{\xi^2+\eta^2}$ near $\rightsonic$, and $r=\sqrt{(\xi-\leftu)^2+\eta^2}$ near $\leftsonic$.
\item[\rm (v)] The limits $\displaystyle\lim_{(\xi,\eta)\to \lefttop\atop(\xi,\eta)\in\Om}D^2\vphi$ and $\displaystyle\lim_{(\xi,\eta)\to \righttop\atop(\xi,\eta)\in\Om}D^2\vphi$ do not exist.
\end{itemize}
\end{theorem}

\begin{remark}
{\rm Assertion (iii)} of {\rm Theorem \ref{theorem-3}} says that $\varphi$ in $\Omega$
is $C^{2,\alpha}$ up to the sonic arcs $\leftsonic$ and $\rightsonic$ away
from the points  $\lefttop$ and $\righttop$.
{\rm Assertion (ii)} of {\rm Theorem \ref{theorem-3}} combined with {\rm (iii)} of {\rm Problem \ref{fbp}}
imply that the function $\vphi_*$ in \eqref{extsol} is $C^{1,1}$ across
the sonic arcs $\rightsonic$ and $\leftsonic$.
{\rm Assertion (iv)} of {\rm Theorem \ref{theorem-3}}  implies that the function $\vphi_*$
is not $C^2$ across the  sonic arcs $\rightsonic$ and $\leftsonic$
since there is a jump of second derivative of $\varphi$, which implies that the $C^{1,1}$-regularity is optimal.
\end{remark}

\begin{remark}[Weak solutions of Problem \ref{fbp} for $\beta=0$]
For $\beta=0$, $\vphi=\rightvphi$ is the unique weak solution of {\rm Problem \ref{fbp}}
satisfying {\rm properties (i)--(v)} of {\rm Theorem \ref{theorem-1}},
where $\rightvphi$ is defined by \eqref{2-k}.
\end{remark}

The following theorem easily follows from Theorem \ref{theorem-1}, Lemma \ref{lemma-parameters},
and the argument after \eqref{extsol}.

\medskip
\begin{theorem}
\label{theorem-2}
For any given $\gam\ge 1$ and any $(\iu, \theta)\in \mathfrak{P}_1$, {\rm Problem \ref{problem-2}}
has a global weak solution $\vphi_*$ of the structure as in {\rm Figure \ref{fig:global-structure}}
with $\vphi_*$ being continuous in $\Lambda$ and $C^1$ across $\leftsonic\cup\rightsonic$.
\end{theorem}

\medskip
\section{Overview of the proof of Theorem \ref{theorem-1}}\label{section-overview}

The key to establish Theorem \ref{theorem-1} is in the monotonicity properties
of $\ivphi-\vphi$ for a weak solution $\vphi$ of Problem \ref{fbp}.
It implies that the shock is Lipschitz graph in a cone of directions and thus fixes geometry of the problem,
among other consequences. Another key property is following.
Fix $\gam \ge 1$, $\iv>0$, and $\beta_*\in(0,\betac(\iv))$.
Then
there exists a constant $C>0$ depending only on $\gam, \iv$, and $\beta_*$
such that, for any $\beta\in(0,\beta_*]$, a corresponding weak solution $\vphi$
with properties (i)--(v) in Theorem \ref{theorem-1} satisfies
\begin{equation*}
{\rm dist}(\shock, B_1(0,-\iv))\ge \frac 1C>0.
\end{equation*}
This inequality plays an essential role to achieve the ellipticity of equation \eqref{2-1}
in $\Om$. Once the ellipticity is achieved, then we obtain various apriori estimates of $\vphi$ so that
we can employ a degree theory
to establish
the existence of a weak solution for all $\beta\in(0,\beta_*]$.
Since the choice of $\beta_*$ is arbitrary in $(0,\betac(\iv))$,
the existence of a weak solution for any $\beta\in(0,\betac(\iv))$ can be established.

\smallskip
For the rest of this paper, we outline the proof of Theorem \ref{theorem-1}.

\subsection{Monotonicity Property and Its Consequences}
\label{monotSection}
Given $\iv>0$ and $\beta\in(0,\beta_*]$,  set $\leftvec:=(\cos\beta, \sin\beta)$.
Then $\leftvec$ is a unit tangent to $\leftshock$ and satisfies $\leftvec\cdot \hat{\xi}>0$.

\begin{lemma}
\label{lemma-1}
For a fixed $(\iv, \beta)\in \mathfrak{P}_2$, if $\vphi$ is a weak solution
with {\rm properties (i)--(v)} in {\rm Theorem \ref{theorem-1}}, then it satisfies
\begin{align*}
& \der_{\leftvec}(\ivphi-\vphi)< 0\;\;\text{in}\;\;\ol{\Om}\setminus\ol{\leftsonic},\qquad
\der_{\xi}(\ivphi-\vphi)>0\;\;\text{in}\;\;\ol{\Om}\setminus \ol{\rightsonic},
\\
&\der_{\xi}(\vphi-\rightvphi)\le 0, \qquad \der_{\eta}(\vphi-\rightvphi)\le 0\;\;\text{in}\;\;\ol{\Om},\\
&\der_{\leftvec}(\vphi-\leftvphi)\ge 0,\qquad \der_{\eta}(\vphi-\leftvphi)\le 0\;\;\text{in}\;\;\ol{\Om}.
\end{align*}
\end{lemma}

For $\beta\in(0,\beta_*]$, and two unit vectors $\leftvec$ and $\rightvec$ in $\R^2$, define
\begin{equation*}
\mathrm{Cone}(\leftvec, \rightvec):=\{a_1\leftvec
+a_2(-\hat{\xi})\,:\, a_1, a_2\ge 0\}.
\end{equation*}
Note that the vectors $\leftvec$ and $\rightvec$ are not parallel if  $\beta\in(0,\beta_*]$,
thus, $\mathrm{Cone}(\leftvec, \rightvec)$ has non-empty interior.
Let $\mathrm{Cone}^0(\leftvec, \rightvec)$ be the interior of  $\mathrm{Cone}(\leftvec, \rightvec)$.

\begin{remark}
For $\beta=0$, we have  $\leftvec=\hat{\xi}$. Then, for  $\beta=0$,
we define $\mathrm{Cone}(\leftvec, \rightvec)$ by the upper half-plane, that is,
\begin{equation*}
\mathrm{Cone}(\leftvec, \rightvec):=\{\bm u\in \R^2\,:\, \bm u\cdot \hat{\eta}\ge 0\},
\end{equation*}
and let $\cone$ be the interior of the upper-half plane $\{\bm u\in \R^2\,:\, \bm u\cdot\hat{\eta}>0\}$.
Note that this is consistent with the definition for $\beta>0$ in the sense that
$(\mathrm{Cone}^0(\leftvec, \rightvec))_{|\beta} \to (\mathrm{Cone}^0(\leftvec(\beta=0), \rightvec))_{|\beta=0}$
as $\beta\to 0+$, where the convergence is locally in the Hausdorff metric.
\end{remark}

Hereafter, we assume that $\vphi$ is a weak solution satisfying
properties (i)--(v) in Theorem \ref{theorem-1}, unless otherwise specified.
By Lemma \ref{lemma-1}, we have
\begin{equation}
\label{3-c1}
\begin{split}
&\der_{\bm e}(\ivphi-\vphi)<0\qquad \text{in}\;\;\ol{\Om}\;\;\text{for all}\;\;\bm e\in \cone,\\
&\der_{\bm e}(\vphi-\rightvphi)\ge 0\qquad \text{in}\;\;\ol{\Om}
 \;\;\text{for all}\;\;\bm e\in\mathrm{Cone}(-\hat{\xi}, -\hat{\eta}),\\
&\der_{\bm e}(\vphi-\leftvphi)\ge 0\qquad \text{in}\;\;\ol{\Om}
\;\;\text{for all}\;\;\bm e\in\mathrm{Cone}(\leftvec, -\hat{\eta}).
\end{split}
\end{equation}

By \eqref{1-i} and the first inequality in \eqref{3-c1}, if $\vphi$ is a weak solution
with properties (i)--(v) in Theorem \ref{theorem-1} for $(\iv, \beta)\in \mathfrak{P}_2$,
there exists a function $\eta=\fshock(\xi)$ satisfying
\begin{itemize}
\item[(i)] $\shock=\{(\xi,\eta):\xi\in(\xi_{\lefttop}, \xi_{\righttop}), \eta=\fshock(\xi)\}$,
where $\xi_{P_j}$ is the $\xi$--coordinate of the point $P_j$ for $j=1,2$;
\item[(ii)] $\fshock(\xi)$ is strictly monotone in $\xi\in(\oxi, \nxi)$:
\begin{equation}
\label{3-c5a}
\fshock'(\xi)> 0\qquad\text{for all}\;\;\xi\in(\oxi, \nxi).
\end{equation}
\item[(iii)]
there exists a constant $C_1$ depending only
on $\gam, \iv$, and $\beta_*$ such that
\begin{equation}
\label{8-15}
\sup_{\oxi<\xi<\nxi}|f'_{shock}(\xi)|\le C_1
\end{equation}
\end{itemize}
From \eqref{3-c5a}, it easily follows that
\begin{equation*}
\inf_{\beta\in[0,\beta_*]}{\rm dist}(\shock, \Wedge)\ge \inf_{[0,\beta_*]}\eta_{\lefttop}>0.
\end{equation*}
Moreover, \eqref{3-c5a} implies that the region $\Om$ enclosed
by $\shock$, $\leftsonic$, $\rightsonic$, and $\Wedge$ is uniformly bounded
for all $\beta\in (0,\beta_*]$. In fact, we have
\begin{equation*}
\Om\subset\{(\xi,\eta):\leftu-\leftc<\xi<\rightc, 0<\eta<\neta\}.
\end{equation*}
From this, we obtain a constant $C>0$ depending only on $\gam, \iv$, and $\beta_*$
such that, if $\vphi$ is a weak solution with properties (i)--(v) in Theorem \ref{theorem-1}
for $\beta\in (0,\beta_*]$, then \begin{equation}
\label{3-c5}
\begin{split}
& \Om\subset B_C(\bm 0),\\
&\sup_{\Om}|\vphi|\le C,\qquad \|\vphi\|_{C^{0,1}(\ol{\Om})}\le C,\\
&(\frac{2}{\gam+1})^{\frac{1}{\gam-1}}\le \rho \le C\;\;\text{in}\;\;\Om,
\qquad 1 \le \rho\le C\;\;\text{on}\;\;\shock.
\end{split}
\end{equation}

The importance of the last two inequalities in \eqref{3-c1} will
be mentioned in Section \ref{UnifEstSection}.

\subsection{Uniform lower bound for the distance between the shock and the sonic circle of the upstream state}\label{UnifLowerBdShockSonicSection}

We can establish the following crucial proposition.

\begin{proposition}
\label{proposition-4}
Given $\gam\ge 1$ and $\iv>0$,
if $\vphi$ is a weak solution for $\beta\in(0,\beta_*]$
with {\rm properties (i)--(v)} in {\rm Theorem \ref{theorem-1}}, then
\begin{align}
\label{new-dec-2}
&{\rm dist}(\shock, B_{1}(0,-\iv))\ge \frac 1{C}>0.
\end{align}
\end{proposition}

Proposition \ref{proposition-4} is essential to obtain uniform ellipticity constants
of equation \eqref{2-1}.
More precisely, we first write equation \eqref{2-1} with $\rho$ given by \eqref{1-p} and $\irho=1$ as
\begin{equation}
\label{8-48}
{\rm div} \mcl{A}(D\vphi,\vphi,\xi,\eta)+\mcl{B}(D\vphi,\vphi,\xi,\eta)=0
\end{equation}
for ${\bf{p}}=(p_1, p_2)\in \R^2, z\in \R$, and $(\xi,\eta)\in\R^2$, where
\begin{equation}
\label{8-49}
\mcl{A}({\bf p}, z,\xi,\eta)\equiv \mcl{A}({\bf p},z):=\rho(|{\bf p}|^2,z){\bf p},
\quad \mcl{B}({\bf p},z ,\xi,\eta)\equiv \mcl{B}({\bf p},z):=2\rho(|{\bf p}|^2,z).
\end{equation}
Then, using \eqref{new-dec-2}, we can find a constant $C$ such that,
if $\vphi$ is a weak solution for $\beta\in(0,\beta_*]$ and
a set $U\subset \ol{\Om}$ satisfies ${\rm dist}(U, \leftsonic\cup\rightsonic)\ge d>0$, we have
\begin{equation}
\label{9-32}
\frac dC|\kappa|^2\le \sum_{i,j=1}^2\mcl{A}_{p_j}^i\bigl(D\vphi(\xi, \eta),
\vphi(\xi, \eta)\bigr)\kappa_i\kappa_j\le C|\kappa|^2
\end{equation}
for $(\xi, \eta)\in U$ and any $\bm{\kappa}=(\kappa_1, \kappa_2)\in \R^2$.
The important part is that such a constant $C$ can be chosen depending
only on $\gam, \iv$, and $\beta_*$, but independent of $\vphi$ and $d$;
therefore we obtain the uniform estimates of weak solutions of
Problem \ref{fbp} for $\beta\in(0,\beta_*]$.

\subsection{Uniform Estimates of Global Weak Solutions}\label{UnifEstSection}

We fix $\iv>0, \gam\ge 1$, and $\beta_*\in(0,\betac(\iv))$. Thanks to \eqref{9-32}, we can achieve uniform estimates
of weak solutions $\vphi$ in $\Om$ for $\beta\in(0,\beta_*]$.
Because of the degeneracy of the ellipticity of equation \eqref{2-1}
on $\leftsonic\cup\rightsonic$, we need to consider two cases:

\medskip
Case 1:  Away from $\leftsonic\cup\rightsonic$, where equation \eqref{2-1} is uniformly elliptic;

\smallskip
Case 2:  Near $\leftsonic\cup\rightsonic$, where the ellipticity degenerates.

\subsubsection{$C^{2,\alp}$--estimates away from  $\leftsonic\cup\rightsonic$}
By \eqref{9-32}, equation \eqref{2-1} is uniformly elliptic away from the sonic
arcs  $\leftsonic\cup\rightsonic$.
Then we can show that, for any $\alp\in(0,1)$ and $r>0$,
there exists a constant $C>0$ depending only on $\gam, \iv, \beta_*, \alp$, and $r$ such that,
if $\vphi$ is a weak solution with properties (i)--(v) in Theorem \ref{theorem-1}
for $\beta\in(0,\beta_*]$, then we have
\begin{itemize}
\item[(i)] for any $B_{2r}(P)\subset \Om$,
\begin{equation*}
\|\vphi\|_{C^{2,\alp}(\ol{B_r(P)})}\le C;
\end{equation*}
\item[(ii)] for $P\in \Wedge$ and $(B_{2r}(P)\cap \Om)\cap(\leftsonic\cup\rightsonic)=\emptyset$,
\begin{equation*}
\|\vphi\|_{C^{2,\alp}(\ol{B_r(P)\cap \Om})}\le C.
\end{equation*}
\end{itemize}
For uniform estimates of shock curves $\shock$ and weak solutions $\vphi$ on $\shock$,
we use a partial hodograph transform.
For any $\beta\in(0,\betac)$, we have
$$
\hat{\eta}\in \cone.
$$
Then, by the first inequality in \eqref{3-c1}, we can find constants $\delta$ and $\sigma>0$
depending only on $\gam, \iv$, and $\beta_*$ such that, if $\vphi$ is a weak solution
for $\beta\in(0,\beta_*]$, then
\begin{equation}
\label{5-a1}
\der_{\eta}(\ivphi-\vphi)\le -\delta\qquad \text{in}\;\;\mcl{N}_{\sigma}(\shock)\cap \Om
\end{equation}
for $\mcl{N}_{\sigma}(\shock)=\{(\xi, \eta)\in \R^2: {\rm dist}\bigl((\xi, \eta), \shock\bigr)<\sigma\}$.
We introduce a coordinate system $(\xi', \eta')$ with $\xi'=\xi$ and a function
$v(\xi', \eta')$ such that $v$ satisfies
\begin{equation}
\label{4-1}
v\bigl(\xi, (\ivphi-\vphi)(\xi, \eta)\bigr)=\eta.
 \end{equation}
By \eqref{5-a1}, such a coordinate system $(\xi', \eta')$ is well defined.
In the new coordinates, we use equation \eqref{2-1} in the set $\mcl{N}_{\sigma}(\shock)\cap \Om$
and the Rankine-Hugoniot condition:
\begin{equation}
\label{r-h}
\rho D\vphi\cdot\nnu_s=D\ivphi\cdot\nnu_s \qquad \mbox{on $\shock$}
\end{equation}
to obtain an equation and a boundary condition for $v$.
Then we obtain an elliptic equation and an oblique boundary condition for
$v$ in the $(\xi', \eta')$--coordinates.
By \eqref{1-i}, $\shock$ becomes a fixed flat boundary; then we obtain
the uniform estimates of $v$ for $\beta\in[0,\beta_*]$.
By \eqref{4-1}, we obtain $\shock=\{\eta=v(\xi, 0)\}$;
thus the uniform estimates of $v$ imply that, for any $k\in \mathbb{N}$ and $d>0$, there are constants $s, C_k(d)$, and $\hat C_k(d)$
depending on $\gam, \iv,\beta_*$, and $d$ such that,
if $P(\xi_P,\eta_P)\in \shock$ and ${\rm dist}(P,\leftsonic\cup\rightsonic)\ge d$, then
\begin{equation}
\label{10-a}
\begin{split}
&|D^kf_{shock}(\xi_P)|\le C_k(d),\qquad\,\,
|D^k_{(\xi,\eta)}\vphi|\le \hat C_k(d)\;\;\;\text{on}\;\;B_s(P)\cap \Om,
\end{split}
\end{equation}
where $C_k(d)$ and $\hat{C}_k(d)$ depend additionally on $k$.

\subsubsection{Weighted $C^{2,\alp}$--estimates near $\leftsonic\cup\rightsonic$}
\label{subsubsec-sonic-est}
Near the sonic arcs $\leftsonic\cup\rightsonic$,
it is convenient to use the coordinates in which the sonic arcs are flattened.
For that reason, we introduce the $(x,y)$--coordinates as follows:
Choose $\eps>0$ small.
\begin{itemize}
\item For $(\xi, \eta)\in\mcl{N}_{\eps}(\rightsonic):=\{(\xi, \eta): {\rm dist}((\xi, \eta),\rightsonic)<\eps\}$,
we set
$$
x:=\rightc-r,\quad  y:=\theta
$$
for the polar coordinates $(r,\theta)$ centered at $(0,0)$;
\item For $(\xi, \eta)\in \mcl{N}_{\eps}(\leftsonic):=\{(\xi, \eta): {\rm dist}((\xi, \eta),\leftsonic)<\eps\}$, we set $$x:=\leftc-r,\quad  y:=\pi-\theta$$ for the polar coordinates $(r,\theta)$ centered at $(\leftu,0)$.
\end{itemize}
From the definition of $(x,y)$ above, it is easy to see that
\begin{equation*}
\begin{split}
&\bigl(\Om\cap\mcl{N}_{\eps}(\Gam)\bigr)\subset\{x>0, y>0\},\qquad \Gam=\ol{\Om}\cap \{x=0\},\\
&\bigl(\{\eta=0\}\cap \mcl{N}_{\eps}(\Gam)\bigr)=\ol{\Om}\cap \{y=0\}
\end{split}
\end{equation*}
for $\Gam=\leftsonic\;\;\text{or}\;\;\rightsonic$.
For  $\leftvphi$ and $\rightvphi$ defined by \eqref{1-21} and \eqref{2-k}, respectively, set
\begin{equation*}
\psi:=\vphi-\max(\leftvphi, \rightvphi).
\end{equation*}
In fact, choosing $\eps$ sufficiently small, we have
\begin{equation*}
\max(\leftvphi, \rightvphi)=
\begin{cases}
\leftvphi&\;\;\text{in}\;\;\mcl{N}_{\eps}(\leftsonic)\cap \Om,\\
\rightvphi&\;\;\text{in}\;\;\mcl{N}_{\eps}(\rightsonic)\cap \Om.
\end{cases}
\end{equation*}
This implies that
\begin{equation*}
\psi=\vphi-\leftvphi\;\;\text{in}\;\;\mcl{N}_{\eps}(\leftsonic)\cap \Om,
\qquad \psi=\vphi-\rightvphi\;\;\text{in}\;\;\mcl{N}_{\eps}(\rightsonic)\cap \Om.
\end{equation*}
Since we seek a weak solution $\vphi$ of Problem \ref{fbp} to be $C^1$
across $\rightsonic\cup\leftsonic$, $\psi$ satisfies
\begin{equation}
\label{5-a4}
\psi(0,y)\equiv0.
\end{equation}
Using the definition of $\psi$, we can rewrite equation \eqref{2-1} as an equation for $\psi$ in
the $(x,y)$--coordinates near $\leftsonic$ or $\rightsonic$ as follows:
\begin{equation*}
A_{11}(D\psi,\psi,x)\psi_{xx}+2A_{12}(D\psi,\psi,x)\psi_{xy}+A_{22}(D\psi,\psi,x)\psi_{yy}+A(D\psi,\psi,x)=0,
\end{equation*}
where $(A_{ij}, A)(\bm p, z, x)$ are smooth with respect to $(\bm p, z, x)$. Then \eqref{9-32} implies
\begin{equation}
\label{5-a3}
\lambda|\bm{\kappa}|^2\le \frac{A_{11}(D\psi,\psi,x)}{x}\kappa_1^2+2\frac{A_{12}(D\psi,\psi,x)}{\sqrt x}\kappa_1\kappa_2+A_{22}(D\psi,\psi,x)\kappa_2^2\le \frac{1}{\lambda}|\bm{\kappa}|^2
\end{equation}
for a constant $\lambda>0$, where $\lambda$ depends only on $\iv, \gam$, and $\beta_*$.
Using the expressions of the coefficients $A_{ij}(\bm p, z, x)$ and the fact that
an ellipticity constant in \eqref{9-32} is proportional to the distance to the sonic
arcs $\leftsonic\cup\rightsonic$, we obtain a constant $\delta>0$ such that,
for any $\beta\in(0,\beta_*]$, a corresponding $\psi$ satisfies
\begin{equation}
\label{6-b5}
\psi_x\le \frac{2-\delta}{1+\gam}x \qquad\text{in}\;\;\mcl{N}_{\eps}(\leftsonic\cup\rightsonic)\cap \Om,
\end{equation}
and such $\delta>0$ depends only on $\gam, \iv$, and $\beta_*$.
Also, by the last two inequalities in \eqref{3-c1} and the definition of the $(x,y)$--coordinates,
one can easily check that $\psi_x\ge 0$ in $\mcl{N}_{\eps}(\leftsonic\cup\rightsonic)\cap \Om$.
Combining this with \eqref{6-b5}, we obtain
\begin{equation}
\label{6-b6}
|\psi_x|\le Cx\qquad\text{in}\;\; \Om\cap \mcl{N}_{\eps}(\rightsonic\cup\leftsonic).
\end{equation}
Combining \eqref{6-b6} with \eqref{5-a4}, we obtain
\begin{equation}
\label{5-a5}
|\psi(x,y)|\le Cx^2\qquad \text{in}\;\;\mcl{N}_{\eps}(\leftsonic\cup\rightsonic)\cap \Om.
\end{equation}
It is important to note that the constant $C$ in \eqref{6-b6} and \eqref{5-a5}
depend only on $\gam, \iv$, and $\beta_*$.
Then the uniform weighted $C^{2,\alp}$--estimates of $\psi$ for $\beta\in(0,\beta_*]$ are
obtained by the local scaling of $\psi$ and the covering argument.

For $P_0\in \mcl{N}_{\eps}(\leftsonic\cup\rightsonic)\cap (\ol{\Om}\setminus (\leftsonic\cup\rightsonic))$,
set $P_0:=(x_0, y_0)$ in the $(x,y)$--coordinates and define
\begin{equation*}
\psi^{(x_0, y_0)}(S,T)=\frac{1}{d^2}\psi(x_0+dS, y_0+\sqrt d T)\qquad \text{with}\;\;d=\frac{x_0}{10}
\end{equation*}
in $Q_1^{(x_0,y_0)}:=\{(S,T)\in(-1,1)^2: (x_0+dS, y_0+\sqrt d T)\in {\Om}\}$.
For the estimates of $\|\psi^{(x_0, y_0)}\|_{C^{2,\alp}(\ol{Q_{1/2}^{(x_0, y_0)}})}$,
we need to consider the case of $P_0\in \shock$ and the case of $P_0\not\in \shock$, separately.
For $P_0\in \shock$, we use the first inequality in \eqref{3-c1} to show that
$$
\der_y(\ivphi-\vphi)\ge \delta \qquad\mbox{near $\leftsonic$ or $\rightsonic$ for some $\delta>0$}.
$$
Then, as in \eqref{4-1}, we can introduce a coordinate system $(x',y')$ and a function $w(x',y')$ satisfying
\begin{equation*}
w(x,(\ivphi-\vphi)(x,y))=y,
\end{equation*}
and we write an equation and a boundary condition for $w$ by using
equation \eqref{2-1} and the boundary condition \eqref{r-h}.
And, we take the same procedure to the function $\ivphi-\max(\leftvphi, \rightvphi)$
to introduce a function $w_0$ and its equation. Then we take the difference
of the equations for $w$ and $w_0$ to obtain an equation for $w-w_0$ and use it to
estimate $w-w_0$.
This procedure will provide the weighted $C^{2,\alp}$--estimates of $\shock$ near $P_0$
and $\|\psi^{(x_0, y_0)}\|_{C^{2,\alp}(\ol{Q_{1/2}^{(x_0, y_0)}})}$.

For  $P_0\in\Om\cup\Wedge $, it is easier to obtain the estimates
of $\|\psi^{(x_0, y_0)}\|_{C^{2,\alp}(\ol{Q_{1/2}^{(x_0, y_0)}})}$.
Then, the supremum of $\|\psi^{(x_0,y_0)}\|_{C^{2,\alp}(\ol{Q^{(x_0,y_0)}})}$
over $(x_0,y_0)\in \mcl{N}_{\eps}(\leftsonic\cup\rightsonic)\cap \Om$
is uniformly bounded for any $\beta\in(0,\beta_*]$. We define
\begin{equation}\label{parabNorms}
 \|\psi\|^{(par)}_{2,\alp,\mcl{N}_{\eps}(\Gam)\cap \Om}
 :=\sup_{(x_0,y_0)\in
 \mcl{N}_{\eps}(\Gam)\cap \Om}\|\psi^{(x_0,y_0)}\|_{C^{2,\alp}(\ol{Q_{1/2}^{(x_0,y_0)}})}
\end{equation}
for $\Gam=\leftsonic\cup\rightsonic$.
We note that the uniform estimate of  $\|\psi\|^{(par)}_{2,\alp,\mcl{N}_{\eps}(\Gam)\cap \Om}$
automatically provides the uniform $C^{1,1}$--estimates of $\psi$
in $\mcl{N}_{\eps}(\Gam)\cap \Om$ for $\beta\in(0,\beta_*]$.

\subsection{Outline of the proof of Theorem \ref{theorem-1}}
\label{subsec-3-4}
To establish Theorem \ref{theorem-1}, that is, the existence of a weak solution of
Problem \ref{fbp} for all $(\iv, \beta)\in \mathfrak{P}_2$,
we apply an iteration procedure. For that purpose, we define an iteration
set $\mcl{K}$ and an iteration map $\mcl{F}$ on $\mcl{K}$ so that a fixed point of $\mcl{F}$ yields a solution of Problem \ref{fbp}.
Then, we apply the Leray-Schauder degree theorem to show the existence of fixed points.

Given $\gam\ge 1$ and $(\iv, \beta)\in \mathfrak{P}_2$,
let $\vphi$ be a solution of Problem \ref{fbp} with a shock $\shock$
satisfying (i)--(v) of Theorem \ref{theorem-1}.
By (ii) of Problem \ref{fbp}, $\shock$ is determined by the solution $\vphi$.
Moreover, by the monotonicity property \eqref{3-c1}, $\shock$ is a graph in $\eta$-direction,
so we can introduce an invertible mapping $\mathfrak{G}$
which maps $\Om$ onto the rectangular domain $\iter:=(-1,1)\times(0,1)$.
Then the function $U=\vphi\circ \mathfrak{G}^{-1}$ is well defined on $\iter$ regardless of $\beta$.
In other words, any solution of Problem \ref{fbp} satisfying (i)--(v)
of Theorem \ref{theorem-1} can be re-defined as a function on $\iter$.
This allows us to perform iteration in a set of functions defined in $\iter$
to prove the existence of a solution of Problem \ref{fbp}.

We define the iteration set $\mcl{K}$ as a product of a parameter set and
a set of functions defined in $\iter$.
Given $\gam\ge 1$ and $\iv>0$, we fix $\beta_*\in (0,\betac(\iv))$ and set
\begin{equation*}
\mcl{K}:=\cup_{\beta\in[0,\beta_*]} \{\beta\}\times \mcl{K}(\beta),
\end{equation*}
where $\mcl{K}(\beta)$ is a set of functions defined in $\iter$.
We define $\mcl{K}(\beta)$ such that, if $u\in \mcl{K}(\beta)$,
there exist a corresponding pseudo-subsonic region $\Om^{(u)}$ and
a function $\vphi^{(u)}$ defined in $\Om^{(u)}$, which satisfies
several properties including the monotonicity properties \eqref{3-c1},
so that it can be a candidate of solution of Problem \ref{fbp}
satisfying (i)--(v) of Theorem \ref{theorem-1}.

Once the iteration set $\mcl{K}$ is defined, the iteration map $\mcl{F}$ is defined as follows:
Given $(\beta, u)\in \mcl{K}$, define the corresponding pseudo-subsonic domain $\Om^{(u)}$,
enclosed by $\rightsonic, \leftsonic$, and $\Wedge$, and an approximate shock $\shock^{(u)}$,
and set up a boundary value problem for an elliptic equation whose ellipticity degenerates
on $\rightsonic\cup\leftsonic$.
Let $\til{\vphi}$ be the solution of the boundary value problem in $\Om^{(u)}$.
The iteration set $\mcl{K}$ will be defined so that such $\til{\vphi}$ exists.
Then we can define a function $\til u$ in $\iter$ from $\til{\vphi}$.
Then the iteration map $\mcl{F}$ is defined by $\mcl{F}(\beta, u)=\til u$.
The boundary value problem for $\til{\vphi}$ is set up so that,
if $u=\til u$, then $\til{\vphi}$ satisfies equation \eqref{2-1}, with an ellipticity
cutoff in a small neighborhood of $\leftsonic\cup\rightsonic$, the Rankine-Hugoniot
conditions \eqref{1-h}--\eqref{1-i} on the shock $\shock$, and the boundary conditions
stated in (iii)--(iv) of Problem \ref{fbp}.
More specifically, note that two conditions in (iii) of Problem \ref{fbp}
are specified on the sonic arcs. Since the sonic arcs are fixed boundaries, this
looks like an overdeterminacy.
Thus, for the iteration problem, we use only the Dirichlet condition
$$
\til{\vphi}=\vphi_{\beta}=\max(\leftvphi, \rightvphi) \qquad\mbox{on $\leftsonic\cup\rightsonic$},
$$
and then we prove that the condition
$D\til{\vphi}=D\vphi_{\beta}$ on  $\leftsonic\cup\rightsonic$ also holds.
In this proof, we use the elliptic degeneracy of
the equation in $\Omega^{(u)}$ near  $\leftsonic\cup\rightsonic$.
Technically, this follows from the estimates of $\tilde\psi=\tilde\vphi-\vphi_{\beta}$
in norms \eqref{parabNorms} near $\leftsonic\cup\rightsonic$.

If $(\beta, u_*)\in \mcl{K}$ satisfies $\mcl{F}(\beta, u_*)=u_*$,
then we call $u_*$ a fixed point.
In order for a fixed point $u_*$ to provide a solution of Problem \ref{fbp},
we need to show that $\vphi_*$ satisfies equation \eqref{2-1}
without the ellipticity cutoff.
Moreover, since we intend to apply the Leray-Schauder degree theorem
to establish Theorem \ref{theorem-1}, we need to prove the following properties:
\begin{itemize}
\item[(i)] If $u_*\in\mcl{K}(\beta)$ is a fixed point of $\mcl{F}$
  for some $\beta\in[0,\beta_*]$, then the corresponding $\vphi_*$ satisfies
   equation \eqref{2-1} without the ellipticity cutoff;
\item[(ii)] For any $\beta\in[0,\beta_*]$, the map $\mcl{F}(\beta, \cdot)$ is compact,
and $\mcl{F}$ is continuous;
\item[(iii)] The iteration set $\mcl{K}$ is open;
\item[(iv)] For any $\beta\in[0,\beta_*]$, there is no fixed point of $\mcl{F}$
on the boundary of the iteration set $\mcl{K}(\beta)$.
\end{itemize}
In proving all the properties above, the apriori estimates in
Sections 4.1 and 4.2
play an essential role. Once (i)--(iv) are verified,
then $\deg(\mcl{F}(\beta, \cdot)-Id, \mcl{K}(\beta),0)$ is a constant
for all $\beta\in[0,\beta_*]$.

To complete the proof of Theorem \ref{theorem-1},
we show that $\deg(\mcl{F}(0, \cdot)-Id, \mcl{K}(0),0)\neq 0$ in two steps.
First, we prove that $\mcl{F}(0, \cdot)$ has the unique fixed point $u^{(normal)}$
in $\mcl{K}(0)$. In fact, we have $u^{(normal)}\equiv 0$.
Then we can conclude that $\deg(\mcl{F}(0, \cdot)-Id, \mcl{K}(0),0)\neq 0$
by showing that $D_u\mcl{F}(0,u^{(normal)})-I$ has the trivial kernel.

\smallskip
As mentioned earlier, once we establish Theorem \ref{theorem-1}, that is,
the existence of global weak solutions of Problem \ref{fbp},
we conclude the existence of a global self-similar weak shock solution
for any admissible wedge angle $\theta_w$ by the one-to-one correspondence between
the parameter sets $\mathfrak{P}_1$ and $\mathfrak{P}_2$.

\smallskip
For the detailed proof, see Bae-Chen-Feldman \cite{BCF2}.

\bigskip
{\bf Acknowledgments.}
The research of
Gui-Qiang G. Chen was supported in part by the National Science
Foundation under Grants
DMS-0935967 and DMS-0807551, the UK EPSRC Science and Innovation
Award to the Oxford Centre for Nonlinear PDE (EP/E035027/1),
the NSFC under a joint project Grant 10728101, and
the Royal Society--Wolfson Research Merit Award (UK).
The work of Mikhail Feldman was
supported in part by the National Science Foundation under Grants
DMS-0800245, DMS-1101260, and the Vilas Award by the University of Wisconsin-Madison.
The research of Myoungjean Bae was supported in part by the initial settlement
research fund provided by POSTECH.

\end{document}